\documentclass[12pt]{article}
\usepackage[margin=1in]{geometry}
\usepackage{amsmath, amssymb, graphicx, enumerate, color, framed, comment, url}
\usepackage{algorithm, algorithmic}
\usepackage[round]{natbib}

\newcommand{\sumu}{\sum_{u \in U}}
\newcommand{\sumv}{\sum_{v \in V}}
\newcommand{\sumi}{\sum_{i \in I}}
\newcommand{\znorm}[1]{\| #1 \|_0}

\newcommand{\ui}{_u^i}
\newcommand{\vi}{_v^i}
\newcommand{\uv}{_{uv}}
\newcommand{\uvi}{_{uv}^i}
\newcommand{\red}[1]{\color{red}#1}
\newcommand{\argmax}{\operatornamewithlimits{argmax}}
\newcommand{\argmin}{\operatornamewithlimits{argmin}}

\title{Sparsity-Constrained Transportation Problem}
\author{
Annie I. Chen\thanks{Department of Electrical Engineering and Computer Science, Massachusetts Institute of Technology, \texttt{anniecia@mit.edu}.} , 
Stephen C. Graves\thanks{Sloan School of Management, Massachusetts Institute of Technology, \texttt{sgraves@mit.edu}.}
}
\date{February 2014}

\begin{document}


\maketitle

\begin{abstract}

We study the solution of a large-scale transportation problem with an additional constraint on the sparsity of inbound flows. Such problems arise in the management of inventory for online retailers that operate with many order fulfillment centers; each stock-keeping unit is typically kept in a limited number of these order fulfillment centers so as to reduce the operational overhead in the supply chain network.

We propose a computationally efficient algorithm that solves this sparsity-constrained optimization problem while bypassing complexities of the conventional integer programming approach. The effectiveness of the algorithm is demonstrated through a series of numerical experiments on synthetic data.


\end{abstract}




\section{Introduction}

There has been a boost in the business of online retail in recent years. The surge in customer demand has led to an expansion of available products online, which in turn again contributes to growth in demand. In addition, due to the unique nature of the virtual platform, products no longer need to be stocked in brick-and-mortar stores to serve customers in the immediate geographic vicinity, but can be distributed across a set of order fulfillment centers and shipped directly on demand. The resulting supply chain networks, being large-scale and integrated, give rise to the need for novel modeling and solution approaches. 

Several important problems in online retail inventory management have been addressed. For example, significant attention has been devoted to efficient fulfillment strategies (i.e., which fulfillment center to ship from when a new order is received) and stocking policies (i.e., how much inventory to stock at each fulfillment center). A review of related topics can be found in \citet{John02}, \citet{Swam03}, and \citet{Agat08}.

The problem of interest in this paper is \emph{inventory positioning}, a decision that precedes considerations for fulfillment and stocking level. Given the demand of products at each geographic zone, the shipping costs between fulfillment centers and demand zones, and the capacity constraints at each fulfillment center, we would like to determine, for each product, which fulfillment centers it should be placed in so as to minimize the overall shipping cost. Without additional constraints, this problem could be modeled as a transporation problem and solved using well-known algorithms for multi-commodity network flow. 
However, due to the large number of distinct stock-keeping units (or items), it is often more practical to limit the number of fulfillment centers that carry each item, so as to reduce the operational overhead for managing the network. We refer to this as the \emph{sparse-inbound constraint}, since it requires the vector of inbound flows in the transportation problem (which correspond to stocking levels at the fulfillment centers) to be sparse, i.e., to have a limited number of nonzero entries. To the best of our knowledge, this is the first systematic investigation of the inventory positioning problem under such conditions. 

Sparsity-constrained optimization problems have been studied in the context of signal processing, particularly for applications where it is necessary to represent high-dimensional signals with sparse estimations that preserve most of the information while requiring much less memory space to store. Several optimization methods have been suggested, including the greedy approaches presented in \citet{Beck12} and \citet{Bahm13}. However, these algorithms are valid only for problems that do not have other constraints besides sparsity, and therefore cannot be directly adapted to our setting where additional network flow constraints need to be satisfied. 

This paper presents an approach to the sparse-inbound transportation problem. The next section provides a mathematical formulation of the problem. The following section introduces a fast algorithm for finding near-optimal solutions. Its effectiveness is then demonstrated with numerical experiments. The final section concludes the work and briefly summarizes our ongoing research.



\section{Problem Formulation}

The sparse-inbound transportation problem can be formulated as a multi-commodity network flow optimization problem on a bipartite graph $G = (U \cup V, E)$. The nodes in the network are partitioned into disjoint sets $U$ and $V$, which represent the order fulfillment centers and the demand zones, respectively. Every edge in $E$ is directed from a node in $U$ to a node in $V$. Let $I$ represent the set of commodities (or items) that flow across the network. We assume that the following information is given:
\begin{itemize}
\item Cost $c\uv$: the cost for sending a unit from fulfillment center $u \in U$ to demand zone $v \in V$.
\item Capacity $l_u$: the maximum total flow (the total amount of items) that each fulfillment center $u \in U$ is able to process. 
\item Demand $z\vi$: the required outbound flow for each item $i \in I$ at each demand zone $v \in V$. 
\item Sparsity $s_i$: the maximum number of fulfillment centers that each item $i \in I$ is allowed to flow through. 
\end{itemize}

Our goal is to find a set of cost-minimizing flows $x = \{x\uvi\}$ and $y = \{y\ui\}$, where $x\uvi$ denotes the flow of item $i$ across edge $(u,v)$, and $y\ui$ is the inbound flow of item $i$ at node $u$. In addition to the standard network flow constraints of flow conservation, edge capacity, and nonnegativity of flows, each inbound flow vector $y^i = (y\ui)_{u \in U}$ also has to satisfy the given sparsity constraint of not having more than $s_i$ nonzero entries. Mathematically, this optimization problem can be described as follows: 

\begin{align}
\mbox{\bf (P)} \quad
\min_{x,y} \quad & \sum_{(u,v) \in E} c\uv \sumi x\uvi \nonumber \\
\mbox{subject to} \quad 
& \znorm{y^i} \leq s_i 		&& \forall i \in I 						&& \mbox{ (inbound flow sparsity) } \label{sparsity}\\
& \sumi y\ui \leq l_u		&& \forall u \in U						&& \mbox{ (inbound capacity) } \label{capacity}\\
& y\ui = \sumv x\uvi 		&& \forall u \in U, \forall i \in I  	&& \mbox{ (inbound flow conservation) } \label{conservation_in}\\
& z\vi = \sumu x\uvi 		&& \forall v \in V, \forall i \in I 		&& \mbox{ (outbound flow conservation) } \label{conservation_out}\\
& x\uvi \geq 0				&& \forall (u,v) \in E, \forall i \in I	&& \mbox{ (nonnegativity) } \label{nonneg_x}\\
& y\ui \geq 0				&& \forall u \in U, \forall i \in I		&& \mbox{ (nonnegativity) } \label{nonneg_y}
\end{align}

\subsection*{Limitations of the MIP approach}

The complexity of problem {\bf(P)} arises mainly from the sparsity constraints \eqref{sparsity}. The conventional approach is to model such constraints as the following mixed integer program (MIP): 
\begin{align}
\mbox{\bf (P-MIP)} \quad
\min_{x,y} \quad & \sum_{(u,v) \in E} c\uv \sumi x\uvi  \nonumber\\
\mbox{subject to} \quad 
& \red{ y\ui \leq M b\ui 	&& \forall i \in I, u \in U	}	\label{mip1}\\
& \red{ \sumu b\ui \leq s_i 	&& \forall i \in I }				\label{mip2}\\
& \sumi y\ui \leq l_u		&& \forall u \in U						&& \mbox{ (inbound capacity) } \nonumber\\
& y\ui = \sumv x\uvi 		&& \forall u \in U, \forall i \in I  	&& \mbox{ (inbound flow conservation) } \nonumber\\
& z\vi = \sumu x\uvi 		&& \forall v \in V, \forall i \in I 		&& \mbox{ (outbound flow conservation) } \nonumber\\
& x\uvi \geq 0				&& \forall (u,v) \in E, \forall i \in I	&& \mbox{ (nonnegativity) } \nonumber\\
& y\ui \geq 0				&& \forall u \in U, \forall i \in I		&& \mbox{ (nonnegativity) } \nonumber\\
& \red{ b\ui \in \{0,1\} } \nonumber
\end{align}
In this formulation, $b\ui$ are binary variables indicating whether or not item $i$ flows through fulfillment center $u$, and $M$ is a large positive constant (which can be set to, for example, the sum of demands for all items in all demand zones). The constraints \eqref{mip1} ensure that $b\ui = 1$ if $y\ui > 0$, and \eqref{mip2} limits the total number of edges with positive flows.

Standard MIP solvers can be applied to obtain the optimal solution of {\bf(P-MIP)}. However, due to the combinatorial nature of the problem, solving the MIP requires intensive computation and does not scale up to practical applications which may have hundreds of fulfillment centers and demand zones as well as thousands or even millions of items. It is therefore necessary to consider heuristic approaches that can find near-optimal solutions quickly. 


\section{The Algorithm}

We now present an efficient two-stage algorithm that finds near-optimal solutions without using integer variables. We begin by explaining an important concept that frames the algorithm as a search algorithm; then, we present details about the two stages of the algorithm, followed by a discussion on trade-offs captured by the selection of parameters. 

\subsection{Conceptual description}

The proposed algorithm is based on the following key observation about the set of feasible solutions of {\bf(P)}: though not a convex set, the feasible set of {\bf(P)} can be expressed as the union of convex sets (in particular, polyhedrons). Indeed, the sparsity constraint \eqref{sparsity} is equivalent to requiring $N_i := |U|-s_i$ of the variables in the vector $y^i = (y\ui)_{u \in U}$ to be 0: 
\begin{align*}
& \znorm{y^i} \leq s_i \quad \forall i \in I 
\\ \Leftrightarrow \quad
& \left( y_{u_1^i}^i = y_{u_2^i}^i = \cdots = y_{u_{N_i}^i}^i = 0  \quad \forall i \in I \right) 
\mbox{ or } 
\left( y_{(u_1^{i})'}^i = y_{(u_2^{i})'}^i = \cdots = y_{(u_{N_i}^{i})'}^i = 0  \quad \forall i \in I \right) 
\mbox{ or } \cdots 
\\ \Leftrightarrow \quad 
& \bigcup_{ \{\tilde{U}_i, i \in I : |\tilde{U_i}| = N_i \} } \left( y\ui = 0 \quad \forall u \in \tilde{U_i}, i \in I \right). 
\end{align*}

Each statement in parenthesis above can replace \eqref{sparsity} in {\bf(P)} to result in a linear program (with a polyhedral feasible set) whose optimal solution satisfies the original sparsity constraint. We shall henceforth refer to such a linear program as a \emph{sparse-LP}. For any given sparse-LP, we say that the node $u$ is \emph{inactive} for item $i$ if the constraint $y\ui = 0$ is present; otherwise, it is \emph{active} for item $i$. 

As we can see from the expression above, there are $\left( \begin{smallmatrix} |U| \\ N_i \end{smallmatrix} \right)$ ways to choose the set $\tilde{U}_i$ for each item $i$. Therefore, there are $\prod\limits_{i \in I} \left( \begin{smallmatrix} |U| \\ N_i \end{smallmatrix} \right)$ possible sparse-LPs (each corresponding to a parenthesized term above), and the union of their feasible sets is the feasible set of {\bf(P)}. Moreover, the optimal solution of {\bf(P)} is the overall optimal solution among all sparse-LPs. 

Our algorithm is essentially a search algorithm over the optimal solutions of all possible sparse-LPs. Since the number of sparse-LPs increases exponentially with $|I|$, it is impossible to do an exhaustive search; however, as we shall see, our algorithm demonstrates that good heuristics could indeed go a long way. 

\subsection{The sparsify-improve algorithm}

There are two consecutive stages in the algorithm, \emph{sparsify} and \emph{improve}: 

\begin{enumerate}

\item In the \emph{sparsify} stage, a sparse-LP is found by progressively ``deactivating'' nodes in a relaxed problem, i.e., adding constraints of the form $y\ui = 0$ for some chosen $u$ and $i$. 

Specifically, the algorithm starts with the relaxed linear program of {\bf(P)} without the sparsity constraint \eqref{sparsity}. In each successive iteration, it solves for the optimal solution of the current linear program, finds $k_1$ active edges $y\ui$ with the smallest flows, and deactivates each of them by adding the constraint $y\ui = 0$ to the current linear program. This process is repeated until all items are left with no more than $s_i$ active nodes and the sparsity constraint is satisfied. 

Note that in general, due to the capacity constraints, not all sparse-LPs are feasible, and there is no guarantee that the sparsify stage will produce a feasible solution. However, we have found empirically that if {\bf(P)} is feasible and the feasible set is not too small, then this search process will find a feasible sparse-LP with high likelihood.

\item Once a (not necessarily optimal) sparse-LP is found, the algorithm then enters the \emph{improve} stage, in which it seeks to improve the solution by exploring other sparse-LPs that are obtained by swapping inactive nodes in the current sparse-LP with active nodes of the same item. 

In particular, in each iteration, the algorithm checks up to $k_2$ inactive inbound edges, in decreasing order of the dual variable magnitude for the constraint $y\ui = 0$ (or in increasing order of dual variable values, since the values are negative), to see if activating the edge results in an improvement. Dual variables are a natural choice for a heuristic rule, because they reflect the marginal cost of each constraint, i.e., the rate at which the optimal cost would change if the constraint $y\ui = 0$ were relaxed to allow item $i$ to flow across node $u$. 

For each inactive edge selected for checking, the algorithm removes the constraint $y\ui = 0$, thereby activating node $u$ for item $i$ and obtaining a relaxed problem. In order to satisfy the sparsity constraint, however, it is necessary to deactivate another node; we choose deactivate the node for the same item that has the smallest flow in the relaxed problem. The optimal solution is then updated by resolving the new sparse-LP. 

The iteration ends once an improvement is found, or if there are no improvements found after checking $k_2$ inactive edges with the largest duals. In the former case, a new iteration is initiated; in the latter case, the algorithm terminates and returns the sparse solution that has been found. 

\end{enumerate}

We summarize this in pseudocode on the next page. 
\begin{algorithm}
\caption{The Sparsify-Improve Algorithm}
\label{algo}
\begin{algorithmic}
\algsetup{indent=1cm}

\vspace{0.5cm}
	\STATE {\bf Function definitions: }
	\STATE \hspace{1cm} $sol(p) :=$ optimal solution of problem $p$ (a vector of flows) 
	\STATE \hspace{1cm} $val(p) :=$ optimal value of problem $p$ (the minimum cost)
	\STATE \hspace{1cm} $\argmin_{i}[k]\{ y_i : i \in I \} :=$ indices of the $k$ smallest elements

\vspace{0.5cm}
	\STATE {\bf Input: } sparse-inbound transportation problem {\bf(P)}

\vspace{0.5cm}
	\STATE {\bf Initialize: }
	\STATE \hspace{1cm} $lp \leftarrow$ {\bf(P)} without \eqref{sparsity}

\vspace{0.5cm}
	\STATE {\it // Stage 1: Sparsify}
	\WHILE{$lp$ is not sparse}
		\STATE $smallest\_flows \leftarrow \argmin_{(u,i)} [k_1] \{y_u^i : y_u^i \in sol(lp), (u,i) \mbox{ active in } lp \}$
		\FOR{$(u,i)$ in $smallest\_flows$} 
			\STATE add constraint $y_u^i = 0$ to $lp$ \quad // deactivate $(u,i)$ for $lp$
		\ENDFOR
	\ENDWHILE

\vspace{0.5cm}
	\STATE {\it // Stage 2: Improve}
	\STATE $improving \leftarrow \TRUE$
	\WHILE{$improving$} 
		\STATE $best\_duals \leftarrow \argmax_{(u,i)} [k_2] \{ | dual(y_u^i = 0) | : y_u^i \in sol(lp), (u,i) \mbox{ inactive} \}$
		\STATE $new\_LPs \leftarrow \phi$
		\FOR{$(u,i) \in best\_duals$}
			\STATE $lp' \leftarrow lp$
			\STATE remove constraint $y_u^i = 0$ from $lp'$ \quad // activate $(u,i)$ for $lp'$
			\STATE $(u',i') \leftarrow \argmin_{(u,i)} \{ y_u^i : y_u^i \in sol(lp'), (u,i) \mbox{ active}\}$
			\STATE add constraint $y_{u'}^{i'} = 0$ to $lp'$ \quad // deactivate $(u',i')$ for $lp'$
			\IF{$val(lp') < val(lp)$} 
				\STATE $lp \leftarrow lp'$
				\STATE go to the top of the while loop
			\ENDIF
		\ENDFOR
		\STATE $improving \leftarrow \FALSE$ 
	\ENDWHILE
	
\vspace{0.5cm}
	\STATE {\bf Output:} $sol(lp)$
\end{algorithmic}
\end{algorithm}

\subsection{Parameter selection}
There are two parameters for this algorithm: $k_1$, the number of nodes to deactivate in each iteration of the \emph{sparsify} stage, and $k_2$, the maximum number of alternative sparse-LPs to consider in the \emph{improve} stage. 

The choice of $k_1$ determines the number of iterations in the \emph{sparsify} stage. The algorithm begins with a relaxed linear program where all $|U|\times|I|$ nodes and items are active, but only $\sumi s_i$ are to be kept in the final sparse-LP. Deactiving the remaining $N := |U||I| - \sumi s_i$ nodes requires $\frac{N}{k_1}$ iterations, which decreases as $k_1$ increases. However, as $k_1$ grows, the search also becomes ``coarser'', and the resulting sparse-LP would likely be less optimal, since fewer intermediate solutions are considered on the search path. Thus, there is a trade-off between optimality and computation time for the choice of $k_1$. Opting for computation time seems to be the better approach, since the \emph{sparsify} stage is only used for finding an initial sparse-LP, which would then be further refined in the \emph{improve} stage. Therefore, we suggest scaling $k_1$ according to the problem size. 

As for $k_2$, the likelihood of finding an improved solution increases with $k_2$, and therefore the \emph{improve} stage is expected to terminate with fewer iterations. However, the running time of each iteration could also increase linearly with $k_2$, since up to $2k_2$ linear programs are solved in each iteration. In addition, the ``marginal likelihood'' of finding a better solution decreases as $k_2$ grows large--- in other words, loosely speaking, the best improvement possible is likely to have a large dual variable value, so if $k_2$ is ``large enough'' that it already includes the best improvement possible, further increasing it does not help with finding a better solution. Due to these complicated factors, there is no obvious choice for $k_2$, and it should be selected based on the allowable computation time.


\section{Numerical Experiments}

In this section, we demonstrate the effectiveness of our algorithm by presenting numerical experiment results for problems of various sizes. As we shall see, the algorithm significantly outperforms the MIP approach in its capability for handling large-scale problems, finding near-optimal solutions in a much shorter amount of time.  

\subsection{Setup}

In our experiments, the transportation networks are generated by placing nodes (order fulfillment centers and demand zones) uniformly at random on the unit square $[0,1) \times [0,1)$. To ensure feasiblity, an edge is created to connect every warehouse to every demand zone. Characteristics of the network are determined as follows: 
\begin{itemize}
\item The cost $c\uv$ for edge $(u,v)$ is set to be the Euclidean distance between fulfillment center $u$ and demand zone $v$. 
\item The capacity $l_u$ of each fulfillment center $u$ is set to be a large number such that all sparse-LPs are feasible. In our case, all capacities were set to $2 \sumi \frac{\sumv z\vi}{s_i}$, which is twice the total demand averaged across the number allowable fulfillment centers for all items. 
\item The demand $z\vi$ for each item $i$ at each demand zone $v$ is a randomly generated integer between 10 and 1000. 
\item The number of allowable fulfillment centers is set to $s_i=5$ for every item. 
\end{itemize}

The experiment is performed on problems with $|U|=30$ fulfillment centers, $|V| = 100$ demand zones, and a varying numbers of items ranging from 1 to 64 in powers of 2.  
Parameters for the \emph{sparsify} and \emph{improve} stages are set to $k_1 = \lceil\sqrt{|I|}\rceil$ and $k_2 = 20$. 

The program is written in Python and uses the Gurobi Optimizer. Experiments are executed on a machine with a dual-cord 1.60GHz processor. 

\subsection{Results}
For each value of $|I|$, the experiment is repeated on 10 randomly generated networks. For each problem instance, when possible, we also solve for the exact optimal solution using the MIP approach {\bf(P-MIP)}. The results for the Sparsify-Improve algorithm are compared with that for the MIP according to two performance metrics: computation time and optimality. 

\begin{table}[h!]
\caption{Computation time and optimality of Sparsify-Improve and MIP}
\label{table1}
\begin{center}
\begin{tabular}{|c|c|c|c|c|}
\hline
$|I|$ & \multicolumn{3}{c}{Average computation time} & Average optimality (\%) \\
	  & MIP & Sparsify-Improve & S-I/MIP (\%) & \\ \hline
1		&	0.917 &	1.523 &	166.09 &		2.10		\\ \hline
2		&5.46		&4.013	&73.50	&3.35	\\ \hline
4		&17.26	&11.01	&63.78	&4.74	\\ \hline									
8		&90.13	&39.68	&44.03	&5.41	\\ \hline
16		&376.15	 &111.97	&29.77	&6.35	\\ \hline
32		&997.63	&338.17		&33.90	&5.63	\\ \hline
64		& - 	&1058.281	& - 	& - 	\\ \hline
\end{tabular}
\end{center}
\end{table}						

\begin{figure}[h!]
\centering
\includegraphics[width=0.7\textwidth]{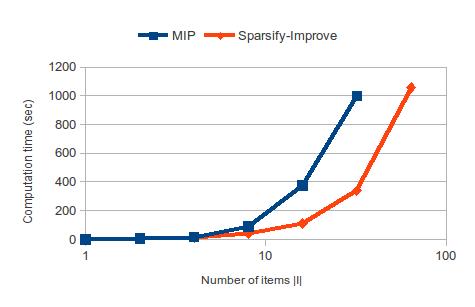}
\caption{Computation time comparison of Sparsify-Improve and MIP}
\label{fig1}
\end{figure}

As we can see from the results summarized in Table \ref{table1} and Figure \ref{fig1}, the Sparsify-Improve algorithm produces solutions that are around 5\% optimal on average, and it takes much less computation time than the MIP (the column labeled ``S-I/MIP'' is the ratio of the average computation time of the Sparsify-Improve algorithm to that of the MIP). Moreover, it also requires less computational resources than the MIP--- for example, in the case where $|I|=64$, the MIP program runs out of memory, and therefore no results are available for comparison.


\section{Conclusion}
In the previous sections, we presented a mathematical formulation of the sparse-inbound transportation problem, which is relevant for online retail inventory positioning problems with sparsity constraints on the number of fulfillment centers. 
We developed a two-stage algorithm for optimizing this model: the \emph{sparsify} stage produces a sparse solution, which is then refined in the \emph{improve} stage. 
The algorithm is based on heuristics, and it significantly outperforms a naive MIP approach, finding a near-optimal solution whose cost is around 5\% optimal with much less computation time. 
Our findings make it possible to efficiently compute sparse inventory positioning plans, thus simplifying the management of large-scale online retail supply chain networks. 

We are currently continuing our investigations on the sparsity-constrained inventory positioning problem. Our pursuits include a better understanding of the optimal solution structure, as well as theoretical analyses of the algorithm. Among other things, this could lead to more concrete performance guarantees, better guidelines for choosing algorithm parameters, and possibly even more effective algorithms for the sparse-inbound transportation problem that could be applied to generic sparsity-constrained optimization problems. We are also exploring the incorporation of additional constraints that reflect other operational considerations, as well as alternative models and algorithms for this problem.



\bibliographystyle{apalike}
\bibliography{bib.bib}

\end{document}